\documentclass[11pt]{article}

\usepackage{amsmath}
\usepackage{amssymb}
\usepackage{epsfig}
\usepackage{amsthm}
\usepackage[mathcal]{eucal}

\newcommand{\F}{\mathcal{F}}
\newcommand{\J}{\mathcal{J}}
\newcommand{\Sf}{\mathcal{S}}

\newcommand{\R}{\mathbb{R}}
\newcommand{\BR}{\bar{\mathbb{R}}}

\newcommand{\inner}[2]{\langle{#1},{#2}\rangle}

\newcommand{\tos}{\rightrightarrows}

\DeclareMathOperator{\clconv}{cl\,conv}

\newtheorem{theorem}{Theorem}[section]

\newtheorem{corollary}[theorem]{Corollary}

\newtheorem{proposition}[theorem]{Proposition}

\title{Maximal monotonicity, conjugation and the duality product in
       non-reflexive Banach spaces}

\author{M. Marques Alves\thanks{IMPA, Estrada Dona Castorina 110, 22460-320
    Rio de Janeiro, Brazil
   ({\tt maicon@impa.br})}\hspace{.5em}\thanks{Partially supported by Brazilian CNPq
    scholarship 140525/2005-0.}
  \and
    B. F. Svaiter\thanks{ IMPA, Estrada Dona Castorina 110, 22460-320 Rio de
    Janeiro, Brazil ({\tt benar@impa.br}) }\hspace{.5em}
    \thanks{Partially supported by CNPq
    grants 300755/2005-8, 475647/2006-8 and by PRONEX-Optimization}
}
\date{}

\begin{document}

\maketitle

\begin{abstract}
  Maximal monotone operators on a Banach space into its dual can be
  represented by convex functions bounded below by the duality
  product.
  It is natural to ask under which conditions a convex function
  represents a maximal monotone operator.
  A satisfactory answer, in the context of reflexive Banach spaces,
  has been obtained some years ago. Recently, a partial result on
  non-reflexive Banach spaces was obtained.
  In this work we study some others conditions which guarantee that a
  convex function represents a maximal monotone operator in
  non-reflexive Banach spaces.
  \\
\end{abstract}

\pagestyle{plain}


\section{Introduction}
Let $X$ be a real Banach space and $X^*$ its topological dual, both
with norms denoted by $\|\cdot\|$. 
The duality product in $X\times X^*$ will be denoted by:
\begin{equation}
  \label{eq:def.pi}
  \pi:X\times X^*\to \R,\quad \pi(x,x^*):=\inner{x}{x^*}=x^*(x). 
\end{equation}
A point to set operator $T:X\tos X^*$ is a relation on $X\times X^*$:
\[ 
 T\subset X\times X^* 
\]
and $T(x)=\{x^*\in X^*\;|\; (x,x^*)\in T\}$. An operator $T:X\tos X^*$
is {\it monotone} if
\[
\inner{x-y}{x^*-y^*}\geq 0,\forall (x,x^*),(y,y^*)\in T
\]
and it is {\it maximal monotone} if it is monotone and maximal (with
respect to the inclusion) in the family of monotone operators of $X$
into $X^*$. 

Fitzpatrick proved constructively that maximal monotone operators are
representable by convex functions.
Before discussing his findigs, let us establish some notation.
We denote the set of extended-real valued functions on $X$ by
$\BR^{X}$.  The {\it epigraph} of $f\in \BR^X$ is defined by 
\[
\mbox{E}(f):=\{(x,\mu)\in X\times \R\,|\,f(x)\leq \mu\}.
\]
  We say that
$f\in \BR^X$ is lower semicontinuous (l.s.c. from now on) if
$\mbox{E}(f)$ is closed in the strong topology of $X\times \R$.

Let $T:X\tos X^*$ be maximal monotone. The \emph{Fitzpatrick
function of $T$} is~\cite{Fitz88}
\begin{equation} \label{FitzIntro}
\varphi_T\in \BR^{X\times X^*},\quad 
 \varphi_{T}(x,x^*):=\sup_{(y,
      y^*)\in T} \inner{x-y}{y^*-x^*}+\inner{x}{x^*}
\end{equation}
and \emph{Fitzpatrick family} associated with $T$ is
\[
  \F_T:=\left\{ h\in \BR^{X\times X^*}
    \left|
      \begin{array}{ll}
        h\mbox{ is convex and l.s.c.}\\
        h(x,x^*)\geq \inner{x}{x^*},\quad \forall (x,x^*)\in X\times X^*\\
        (x,x^*)\in T
        \Rightarrow 
        h(x,x^*) = \inner{x}{x^*}
      \end{array}
    \right.
  \right\}.
\]
In the next theorem we summarize the Fitzpatrick's results:
\begin{theorem}[\mbox{\cite[Theorem 3.10]{Fitz88}}] \label{th:fitz} 
  Let $X$ be a real Banach space and $T:X\tos X^*$ be maximal
  monotone. Then  for any $h\in \F_T$
  \[
  (x,x^*)\in T\iff h(x,x^*)=\inner{x}{x^*}
  \]
  and $\varphi_T$ is the smallest element
  of the  family $\F_T$.
\end{theorem}
\noindent
Fitzpatrick's results described above were rediscovered by
Mart\'inez-Legaz and Th\'era~\cite{LegTheJNCA01}, and Burachik and
Svaiter~\cite{BuSvSet02}.

It seems interesting to study conditions under which a convex function
$h\in \BR^X$ represents a maximal monotone operator, that is,
$h\in\F_T$ for some maximal monotone operator $T$. Our aim is to
extend previous results on this direction.
We will need some auxiliary results and additional notation for this
aim.

The Fenchel-Legendre conjugate of $f\in \BR^{X}$ is
\[
 f^*\in \BR^{X^*},\quad f^*(x^*):=\sup_{x\in X} \inner{x}{x^*}-f(x).
\]
Whenever necessary, we will
identify $X$ with its image under the canonical injection of $X$ into
$X^{**}$.  
Burachik and Svaiter proved that the family $\F_T$ is invariant under
the mapping
\begin{equation}
  \label{eq:def.jconj}
  \J:\BR^{X\times X^*}\to\BR^{X\times X^*}, \;\J\;h(x,x^*):=h^*(x^*,x).
\end{equation}
This means that
if $T:X\tos X^*$ is maximal monotone, then~\cite{BuSvSet02}
\begin{equation}\label{eq:inv.0} 
 \J(\F_T)\subset \F_T.
\end{equation}
In particular, for any $h\in \F_T$  it holds that 
$ h\geq \pi$ , $ \J h\geq \pi$, 
that is,
\[ h(x,x^*)\geq \inner{x}{x^*},\quad h^*(x^*,x)\geq \inner{x}{x^*},
  \qquad \forall (x,x^*)\in X\times X^*.
\]
So, the above conditions are\emph{ necessary} for a convex function
$h$ on $X\times X^*$ to represent a maximal monotone operator.
Burachik and Svaiter proved that these conditions are also
\emph{sufficient}, in a reflexive Banach space, for $h$ to represent a
maximal monotone operator:
\begin{theorem}[\mbox{\cite[Theorem 3.1]{BuSvProc03}}]\label{th:bs}
 Let $h\in \BR^{X\times X^*}$ be proper, convex, l.s.c. and
 \begin{equation} \label{eq:hhs.0}
 h(x,x^*)\geq \inner{x}{x^*},\quad h^*(x^*,x)\geq \inner{x}{x^*},
  \qquad \forall (x,x^*)\in X\times X^*.
  \end{equation}
  If $X$ is reflexive, then
 \[
  T:=\{(x,x^*)\in X\times X^*\;|\;h(x,x^*)= \inner{x}{x^*}\}
 \]
 is maximal monotone and $h, \J h\in \F_T$.
\end{theorem}
\noindent
Marques Alves and Svaiter generalized Theorem~\ref{th:bs}
 to non-reflexive Banach spaces as follows:
\begin{theorem}[\mbox{\cite[Corollary 4.4 ]{MASvJCA08}}]\label{th:as}
 If $h\in \BR^{X\times X^*}$ is convex and
  \begin{equation} \label{eq:hhs.pd}
  \begin{array}{ll}
    h(x,x^*)\geq \inner{x}{x^*}, & \forall (x,x^*)\in X\times X^*,\\[.4em]
     h^*(x^*,x^{**})\geq \inner{x^*}{x^{**}}, &
      \forall (x^*,x^{**})\in X^*\times X^{**}
  \end{array}
\end{equation}
then
 \[
  T:=\{(x,x^*)\in X\times X^*\;|\;h^*(x^*,x)= \inner{x}{x^*}\}
 \]
 is maximal monotone and $ \J h\in \F_T$. Moreover, if $h$ is l.s.c. then
 $h\in \F_T$.
\end{theorem}
Condition~\eqref{eq:hhs.pd} of Theorem~\ref{th:as} enforces the
operator $T$ to be of type (NI)~\cite{MASvJCA08.2} and is not
necessary for maximal monotonicity of $T$ in a non-reflexive Banach
space. Note that 
   while the weaker condition \eqref{eq:hhs.0} of
 Theorem~\ref{th:bs} is still \emph{necessary} in non-reflexive Banach
spaces for the inclusion $h\in\F_T$, where $T$ is a maximal monotone
operator.
The main result of this paper is another generalization of
Theorem~\ref{th:bs} to  non-reflexive Banach spaces which uses
condition~\eqref{eq:hhs.0} instead of \eqref{eq:hhs.pd}. For obtaining
this generalization, we added a regularity assumption on the domain of
$h$.

If $T:X\tos X^*$ is maximal monotone, it is easy to prove that
$\varphi_T$ is minimal in the family  of all convex functions in $X\times X^*$
which majorizes the duality product.
So, it is natural to ask whether the converse also holds, that is:
\begin{quotation}
  \noindent
  Is any minimal element of this family (convex functions which
  majorizes the duality product) a  Fitzpatrick
  function of some maximal monotone operator?
\end{quotation}
To give a partial answer to this question, Mart\'inez-Legaz and Svaiter proved the following results, which we will use latter on:
\begin{theorem}[\mbox{\cite[Theorem 5]{LegSvProc08}}] \label{th:ls}
 Let $\mathcal{H}$ be the family of convex functions in $X\times X^*$ which
 majorizes the duality product:
 \begin{equation}
   \label{eq:def.h}
   \mathcal{H}:=
  \{h\in \BR^{X\times X^*}\,|\,h\; \mbox{is proper, convex and}\; h\geq \pi\}. 
 \end{equation}
 The following statements holds true:
 \begin{enumerate}
 \item The family $\mathcal{H}$ is (donward) inductively ordered;
  \item For any $h\in \mathcal{H}$ there exists a minimal $h_0\in \mathcal{H}$ such that 
        $h\geq h_0$;
      \item Any minimal element $g$ of $\mathcal{H}$ is l.s.c. and
        satisfies $\J g\geq g$.
 \end{enumerate}
\end{theorem}
Note that item 2 is a direct consequence of item 1.
Combining item 3 with Theorem~\ref{th:bs}, Mart\'inez-Legaz and Svaiter
concluded that \emph{in a reflexive} Banach space, any minimal element
of $\mathcal{H}$ is the Fitzpatrick function of some maximal monotone
operator~\cite[Theorem 5]{LegSvProc08}.
We will also present a partial extension of this result for
non-reflexive Banach spaces.

\section{Basic results and notation}

The weak-star topology of $X^*$ will be denoted by
$\omega^*$ and by $s$ we denote the strong topology of $X$.
A function $h\in \BR^{X\times X^*}$  is lower
semicontinuous in the strong$\,\times\,$weak-star topology if
$\mbox{E}(h)$ is a closed subset of $X\times X^*\times \R$ in the
$s\times \omega^*\times |\cdot|$ topology.  

The \emph{indicator function} of $V\subset X$ is $\delta_V$,
$\delta_V(x):=0,\,x\in V$ and $\delta_V(x):=\infty$, otherwise.
The closed convex closure
of $f\in \BR^{X}$ is defined by
\[
 \mbox{cl\,conv}f\in \BR^{X},\quad \mbox{cl\,conv}f(x):=
 \inf\{\mu\in \R\,|\,(x,\mu)\in \mbox{cl\,conv}\,\mbox{E}(f)\} 
\]
where for $U\subset X$, $\mbox{cl\,conv}\,U$ is the closed convex hull
(in the $s$ topology) of $U$.
The {\it effective domain} of a function $f\in \BR^{X}$ is
\[
\mbox{D}(f):=\{x\in X\,|\, f(x)<\infty\},
\]
and $f$ is \emph{proper} if $\mbox{D}(f)\neq
\emptyset$.
If $f$ is proper, convex and l.s.c.,
then $f^*$ is proper.
For $h\in \BR^{X\times X^*}$, we also define
\[
\mbox{Pr}_X\,\mbox{D}(h):=\{x\in X\,|\,\exists\, x^*\in X^*\,|\,(x,x^*)\in \mbox{D}(h)\}.
\]
Let  $T:X\tos X^*$ be  maximal monotone.
In~\cite{BuSvSet02} Burachik and Svaiter defined and studied the
biggest element of $\F_T$, namely, the $\mathcal{S}$-function,
$\mathcal{S}_T\in \F_T$ defined by
\[
\mathcal{S}_T\in \BR^{X\times X^*},\,\quad \mathcal{S}_T:=\sup_{h\in \F_T}\,\{h\},
\]
or, equivalently
\[
  \mathcal{S}_T=\clconv (\pi+\delta_T).
\]
Recall that $\J(\F_T)\subset \F_T$. Additionally~\cite{BuSvSet02}
\begin{equation}\label{eq:inv.1} 
 \J\;\Sf_T=\varphi_T
\end{equation}
and, in a reflexive Banach space, $\J\varphi_T=\mathcal{S}_T$.

In what follows we present the Attouch-Brezis's version of the
Fenchel-Rockafellar duality theorem:

\begin{theorem}[\mbox{\cite[Theorem 1.1]{AttBrez86}}]\label{th:fr}
  Let $Z$ be a Banach space and $\varphi,\psi\in \BR^{Z}$ be proper,
  convex and l.s.c. functions.  If
\begin{equation}\label{eq:cond.fr}
 \bigcup_{\lambda>0}\lambda\,[\emph{D}\,(\varphi)-\emph{D}\,(\psi)],
\end{equation}
is a closed subspace of $Z$, then
\begin{equation}\label{eq:fr}
 \inf_{z\in Z}\varphi(z)+\psi(z)=\max_{z^*\in Z^*}-\varphi^*(z^*)-\psi^*(-z^*).
\end{equation}
\end{theorem}

Given $X,\,Y$ Banach spaces, $\mathcal{L}(Y,X)$ denotes the set of
continuous linear operators of $Y$ into $X$. The range of $A\in
\mathcal{L}(Y,X)$ is denoted by $\mbox{R}(A)$ and the adjoint by
$A^*\in \mathcal{L}(X^*,Y^*)$:
\[
 \inner{Ay}{x^*}=\inner{y}{A^*x^*}\; \forall y\in Y,\,x^*\in X^*,
\]
where $X^*,\,Y^*$ are the dual of $X$ and $Y$, respectively.
The next proposition is a particular case of Theorem 3
of~\cite{Simons-Quadr}. For the sake of completeness, we give the
proof in the Appendix A.

\begin{proposition}\label{prop:2}
 Let $X,\,Y$ Banach spaces and $A\in \mathcal{L}(Y,X)$. For $h\in \BR^{X\times X^*}$,
 proper convex and l.s.c., define $f\in \BR^{Y\times Y^*}$
  \[
   f(y,y^*):=\inf_{x*\in X^*}\,h(Ay,x^*)+\delta_{\{0\}}(y^*-A^*x^*).
  \]
 If
 \begin{equation}\label{eq:cond2}
  \bigcup_{\lambda>0}\lambda\,[\emph{Pr}_{X}\,\emph{D}(h)-\emph{R}(A)],
 \end{equation}
 is a closed subspace of $X$, then
 \[
  f^*(z^*,z)=\min_{u^*\in X^*}\,h^*(u^*,Az)+\delta_{\{0\}}(z^*-A^*u^*).
 \]
\end{proposition}
Mart\'inez-Legaz and Svaiter~\cite{LegSvSet05} defined, for $h\in \BR^{X\times X^*}$ and
$(x_0,x_0^*)\in X\times X^*$, $h_{(x_0,x_0^*)}\in \BR^{X\times X^*}$
\begin{align}\label{eq:def.hx}
 \nonumber
 h_{(x_0,x_0^*)}(x,x^*)&:= h(x+x_0,x^*+x_0^*)-[\inner{x}{x_0^*}+\inner{x_0}{x^*}
    +\inner{x_0}{x_0^*}]\\
                       & = h(x+x_0,x^*+x_0^*)-\inner{x+x_0}{x^*+x_0^*}+\inner{x}{x^*}.
\end{align}
The operation $h\mapsto h_{(x_0,x_0^*)}$ preserves many properties of
$h$, as convexity and lower semicontinuity. Moreover, one can easily prove the 
following Proposition:
\begin{proposition}\label{prop.prh}
Let $h\in \BR^{X\times X^*}$. Then it holds that
\begin{enumerate}
\item $h\geq \pi\iff h_{(x_0,x_0^*)}\geq \pi,\, \forall (x_0,x_0^*)\in X\times X^*$;
\item $\J h_{(x_0,x_0^*)}=(\J h)_{(x_0,x_0^*)}$, \, $\forall (x_0,x_0^*)\in X\times X^*$.
\end{enumerate}
\end{proposition}

\section{Main results}

In the next theorem we generalize Theorem~\ref{th:bs} to non-reflexive
Banach spaces using condition~\eqref{eq:hhs.0} instead of the
condition \eqref{eq:hhs.pd} used in Theorem~\ref{th:as}. For
obtaining this generalization, we added a regularity
assumption~\eqref{eq:c.th:1} on the domain of $h$.
\begin{theorem}\label{th:1}
 Let $h\in \BR^{X\times X^*}$ be proper, convex and
   \begin{equation} \label{eq:hhs.1}
     h(x,x^*)\geq \inner{x}{x^*},\quad h^*(x^*,x)\geq \inner{x}{x^*},
     \qquad \forall (x,x^*)\in X\times X^*.
  \end{equation}
 If
 \begin{equation}\label{eq:c.th:1}
  \bigcup_{\lambda> 0}\lambda\, \emph{Pr}_{X}\,\emph{D}(h),
 \end{equation}
 is a closed subspace of $X$,
 then
 \[
  T:=\{(x,x^*)\in X\times X^*\;|\; h^*(x^*,x)= \inner{x}{x^*}\}
 \]
 is maximal monotone and $\J h\in \F_T$.
\end{theorem}
\begin{proof}
  First, define $\bar h:=\mbox{cl}\,h$ and note that $\bar h$ is
  proper, convex, l.s.c., satisfies~\eqref{eq:hhs.1}, \eqref{eq:c.th:1}
  and $\J \bar h=\J h$. So, it suffices to prove the theorem for the
  case where $h$ is l.s.c., and we assume it from now on in this
  proof. 
  Monotonicity of $T$ follows from Theorem 5 of~\cite{LegSvSet05}.
  Note that for any $x\in X$
  \[
  T(x)=\{x^*\in X^*\;|\; h^*(x^*,x)-\inner{x}{x^*}\leq 0\}.
  \]
  Therefore, $T(x)$ is convex and $\omega^*$-closed.

  To prove maximality of $T$, take $(x_0,x_0^*)\in X\times X^*$ such
  that
\begin{equation}\label{eq:monrel}
 \inner{x-x_0}{x^*-x_0^*}\geq 0,\quad \forall (x,x^*)\in T
\end{equation}
and suppose $x_0^*\notin T(x_0)$. As $T(x_0)$ is convex and
$\omega^*$-closed,
using the geometric version of the Hahn-Banach theorem in 
$X^{*}$ endowed with the $\omega^*$ topology we
conclude that there exists $z_0\in X$ such that
\begin{equation}\label{eq:contr}
 \inner{z_0}{x_0^*}<\inner{z_0}{x^*},\quad \forall\, x^*\in T(x_0).
\end{equation}
Let $Y:=\mbox{span}\{x_0,z_0\}$.
Define $A\in \mathcal{L}(Y,X)$, $A\,y:=y,\,\forall\,y\in Y$ and the convex 
function $f\in \BR^{Y\times Y^*}$,
\begin{equation}\label{eq:def.H.p1}
   f(y,y^*):=\inf_{x*\in X^*}\,h(Ay,x^*)+\delta_{\{0\}}(y^*-A^*x^*).
\end{equation}
Using Proposition~\ref{prop:2} we obtain
\begin{equation}\label{eq:conj.H}
  f^*(y^*,y)=\min_{x*\in X^*}\,h^*(x^*,Ay)+\delta_{\{0\}}(y^*-A^*x^*).
\end{equation}
Using~\eqref{eq:hhs.1},~\eqref{eq:def.H.p1} and~\eqref{eq:conj.H} it is easy to see that
\begin{equation} \label{eq:H}
   f(y,y^*)\geq \inner{y}{y^*},\quad f^*(y^*,y)\geq \inner{y}{y^*},\, \forall (y,y^*)\in Y\times Y^*.
\end{equation}
Define $g:=\J f$. As $Y$ is reflexive we have $\J g=\mbox{cl}\,f$.
Therefore, using~\eqref{eq:H} we also have
\begin{equation} \label{eq:G}
   g(y,y^*)\geq \inner{y}{y^*},\quad g^*(y^*,y)\geq \inner{y}{y^*},\,\forall (y,y^*)\in Y\times Y^*.
\end{equation}
Now, using~\eqref{eq:G} and item 1 of Proposition~\ref{prop.prh} we obtain
\begin{align}\label{eq:1}
\nonumber
g_{(x_0,A^*x_0^*)}(y,y^*)+\frac{1}{2}\|y\|^2+\frac{1}{2}\|y^*\|^2 &\geq
\inner{y}{y^*}+\frac{1}{2}\|y\|^2+\frac{1}{2}\|y^*\|^2\\
     &\geq 0,\quad \forall (y,y^*)\in 
 Y\times Y^*
\end{align}
and
\begin{align}\label{eq:2}
\nonumber
(\J g)_{(x_0,A^*x_0^*)}(y,y^*)+\frac{1}{2}\|y\|^2+\frac{1}{2}\|y^*\|^2 &\geq
\inner{y}{y^*}+\frac{1}{2}\|y\|^2+\frac{1}{2}\|y^*\|^2\\
   &\geq 0,\quad \forall (y,y^*)\in 
 Y\times Y^*.
\end{align}
Using Theorem~\ref{th:fr} and item 2 of Proposition~\ref{prop.prh} we conclude that there exists $(\tilde z,\tilde z^*)\in Y\times Y^*$ 
such that
\begin{equation}\label{eq:3}
 \inf\,g_{(x_0,A^*x_0^*)}(y,y^*)+\frac{1}{2}\|y\|^2+\frac{1}{2}\|y^*\|^2
 +(\J g)_{(x_0,A^*x_0^*)}(\tilde z,\tilde z^*)+\frac{1}{2}\|\tilde z\|^2
  +\frac{1}{2}\|\tilde z^*\|^2=0.
\end{equation}
From~\eqref{eq:1},\eqref{eq:2} and \eqref{eq:3} we have
\begin{equation}\label{eq:4}
 \inf_{(y,y^*)\in Y\times Y^*}\,g_{(x_0,A^*x_0^*)}(y,y^*)+\frac{1}{2}\|y\|^2+\frac{1}{2}\|y^*\|^2
 =0.
\end{equation}
As $Y$ is reflexive, from~\eqref{eq:def.hx},\eqref{eq:4} we conclude that there exists 
$(\hat y, \hat y^*)\in Y\times Y^*$ such that
\begin{equation}\label{eq:5}
 g\,(\hat y+x_0,\hat y^*+A^*x_0^*)-\inner{\hat y+x_0}{\hat y^*+A^*x_0^*}
 +\inner{\hat y}{\hat y^*}+\frac{1}{2}\|\hat y\|^2+\frac{1}{2}\|\hat y^*\|^2
 =0.
\end{equation}
Using~\eqref{eq:5} and the first inequality of~\eqref{eq:G} (and the definition of $g$) we have
\begin{equation}\label{eq:6}
 f^*(\hat y^*+A^*x_0^*,\hat y+x_0)=\inner{\hat y+x_0}{\hat y^*+A^*x_0^*}
\end{equation}
and
\begin{equation}\label{eq:7}
 \inner{\hat y}{\hat y^*}+\frac{1}{2}\|\hat y\|^2+\frac{1}{2}\|\hat y^*\|^2
 =0.
\end{equation}
Using~\eqref{eq:conj.H} we have that there exists $w_0^*\in X^*$ such that
\begin{equation}\label{eq:8}
  f^*(\hat y^*+A^*x_0^*,\hat y+x_0)= h^*(w_0^*,A\,(\hat y+x_0)),\quad \hat y^*+A^*x_0^*=A^*w_0^*.
\end{equation}
So, combining~\eqref{eq:6} and~\eqref{eq:8} we have
\[
 h^*(w_0^*,A\,(\hat y+x_0))=\inner{\hat y+x_0}{A^*w_0^*}=\inner{A(\hat y+x_0)}{w_0^*}.
\]
In particular, $w_0^*\in T(A(\hat y+x_0))$. As $x_0\in Y$, we can use~\eqref{eq:monrel} and the second equality of~\eqref{eq:8} to conclude that
\begin{equation}\label{eq:9}
 \inner{A(\hat y+x_0)-x_0}{w_0^*-x_0^*}=\inner{\hat y}{A^*(w_0^*-x_0^*)}=
 \inner{\hat y}{\hat y^*}\geq 0.
\end{equation}
Using~\eqref{eq:7} and~\eqref{eq:9} we conclude that $\hat y=0$ and $\hat y^*=0$. Therefore,
\[
 w_0^*\in T(x_0),\quad A^*x_0^*=A^*w_0^*.
\]
As $z_0\in Y$, we have $z_0=A\,z_0$ and so
\[
 \inner{z_0}{x_0^*}=\inner{A\,z_0}{x_0^*}=\inner{z_0}{A^*x_0^*}=\inner{z_0}{A^*w_0^*}
 =\inner{A\,z_0}{w_0^*}=\inner{z_0}{w_0^*},
\]
that is,
\[
 \inner{z_0}{x_0^*}=\inner{z_0}{w_0^*},\quad w_0^*\in T(x_0)
\]
which contradicts~\eqref{eq:contr}. Therefore, $(x_0,x_0^*)\in T$ and
so $T$ is maximal monotone and $\J h\in \F_T$.
\end{proof}

Observe that if $h$ is convex, proper and l.s.c. in the
strong$\,\times\,$weak-star topology, then $\J^2h=h$. Therefore, using this
observation we have the following corollary of Theorem~\ref{th:1}:

\begin{corollary}\label{cr:1}
 Let $h\in \BR^{X\times X^*}$ be proper, convex, l.s.c. in the
 strong$\,\times\,$weak-star topology and
  \[
     h(x,x^*)\geq \inner{x}{x^*},\quad h^*(x^*,x)\geq \inner{x}{x^*},
     \qquad \forall (x,x^*)\in X\times X^*.
  \]
 If
 \[
   \bigcup_{\lambda>0}\lambda\, \emph{Pr}_{X}\,\emph{D}(h),
 \]
 is a closed subspace of $X$,
 then
 \[
  T:=\{(x,x^*)\in X\times X^*\;|\;h(x,x^*)= \inner{x}{x^*}\}
 \]
 is maximal monotone and $h,\J h\in \F_T$.
\end{corollary}
\begin{proof}
  Using Theorem \ref{th:1} we conclude that the set
  \[
  S:=\{(x,x^*)\in X\times X^*\;|\; h^*(x^*,x)=\inner{x}{x^*}\}
  \]
  is maximal monotone. Take $(x,x^*)\in S$. As $\pi$ is Gateaux
  differentiable, bounds below $h$ and coincides with $h$ at
  $(x,x^*)$, we have (see Lemma 4.1 of~\cite{MASvJCA08})
  \[
  D\pi(x,x^*)\in \partial\, \J h(x,x^*),
  \]
  where $D\pi$ stands for the Gateaux derivative of $\pi$. As $D\pi(x,x^*)=(x^*,x)$, we conclude that
  \[
  \J h(x,x^*)+\J^2 h(x,x^*)=\inner{(x,x^*)}{(x^*,x)}.
  \]
  Substituting $\J h(x,x^*)$ by $\inner{x}{x^*}$ in the above equation we conclude 
  that $\J^2h(x,x^*)=\inner{x}{x^*}$. Therefore, as $\J^2 h(x,x^*)=h(x,x^*)$,
  \[ 
   S \subset T.
  \]
  To end the proof use the maximal monotonicity of $S$ (Theorem~\ref{th:1})
  and the monotonicity of $T$ (see Theorem 5 of~\cite{LegSvSet05}) to conclude that $S=T$.
\end{proof}

It is natural to ask whether we can drop lower-semicontinuity assumptions. In the context of non-reflexive Banach spaces, we should use
the l.s.c closure in the strong$\,\times\,$weak-star topology. Unfortunately, as
the duality product is not continuous in this topology, it is not clear
whether the below implication holds:
\[
 h\geq\pi\stackrel{?}\Rightarrow \mbox{cl}_{s\times \omega^*}\, h_{}\geq \pi.
\]

\begin{corollary}\label{cor:3}
 Let $h\in\BR^{X\times X^*}$ be proper, convex and
   \[
     h(x,x^*)\geq \inner{x}{x^*},\quad h^*(x^*,x)\geq \inner{x}{x^*},
     \qquad \forall (x,x^*)\in X\times X^*.
   \]
 If
 \[
  \bigcup_{\lambda>0}\lambda\,\emph{Pr}_X\emph{D}(h)
 \]
 is a closed subspace of $X$,
 then
 \[
  \emph{cl}_{s\times \omega^*}\,h\in\F_T,
 \]
 where $\emph{cl}_{s\times \omega^*}$ denotes the l.s.c. closure in the
 strong$\,\times\,$weak-star topology and $T$ is the maximal monotone
 operator defined as in Theorem~\ref{th:1}:
 \[  
 T:=\{(x,x^*)\in X\times X^*\;|\;h^*(x^*,x)=\inner{x}{x^*}\}.
 \]
 In particular, $\emph{cl}_{s\times \omega^*}\,h\geq\pi$.
\end{corollary}
\begin{proof}
  First use Theorem~\ref{th:1} to conclude that $T$ is maximal monotone and
  $\J h\in\F_T$. In particular, 
  \[
  \mathcal{S}_T\geq \J h\geq \varphi_T.
  \]
  Therefore,
  \[
  \J\varphi_T\geq \J^2 h\geq\J\mathcal{S}_T.
  \]
  As $\J\mathcal{S}_T=\varphi_T\in \F_T$ and $\J\varphi_T\in\F_T$, we
  conclude that $\mbox{cl}_{s\times \omega^*}\,h=\J^2h\in\F_T$.
\end{proof}

In the next corollary we give a partial answer for an open
question proposed by Mart\'inez-Legaz and Svaiter in~\cite{LegSvProc08}, in the
context of non-reflexive Banach spaces.
\begin{corollary}\label{cr:2}
Let $\mathcal{H}$ be the family of convex functions on $X\times X^*$
bounded below by the duality product, as defined in~\eqref{eq:def.h}.
If $g$ is a minimal element of $\mathcal{H}$ and
 \[
  \bigcup_{\lambda>0}\lambda\,\emph{Pr}_X\,\emph{D}(g)
 \]
is a closed subspace of $X$,
then there exists a maximal monotone operator $T$ such that $g=\varphi_T$,
where $\varphi_T$ is the Fitzpatrick function of $T$.
\end{corollary}
\begin{proof}
  Using item 3 of Theorem~\ref{th:ls} and Theorem~\ref{th:1} we have that
  \[
   T:=\{(x,x^*)\in X\times X^*\,|\, g^*(x^*,x)=\inner{x}{x^*}\}
  \]
  is maximal monotone, $\J g\in \F_{T}$ and 
  \[ 
    T\subset \{(x,x^*)\in X\times X^*\;|\; g(x,x^*)=\inner{x}{x^*}\}.
  \]
  As $g$ is convex and bounded below by the duality product, using Theorem 5 of~\cite{LegSvSet05},
  we conclude that the rightmost set on the above inclusion is
  monotone. Since $T$ is maximal monotone, the above inclusion holds
  as an equality and, being l.s.c., $g\in\F_T$. To end the proof, note that
  $g\geq\varphi_T\in\mathcal{H}$.
\end{proof}

\appendix

\section{Proof of Proposition~\ref{prop:2}} \label{ap:a}

\begin{proof}[Proof of Proposition~\ref{prop:2}]
Using the Fenchel-Young inequality we have, for any $(y,y^*),(z,z^*)\in Y\times Y^*$ and
$x^*,u^*\in X^*$,
\[
 h(Ay,x^*)+\delta_{\{0\}}(y^*-A^*x^*)+h^*(u^*,Az)+\delta_{\{0\}}(z^*-A^*u^*)
\geq \inner{Ay}{u^*}+\inner{Az}{x^*}.
\]
Taking the infimum over $x^*,u^*\in X^*$ on the above inequality we get
\begin{align}
 \nonumber
 f(y,y^*)+\inf_{u^*\in X^*} h^*(u^*,Az)+\delta_{\{0\}}(z^*-A^*u^*)&\geq
                         \inner{y}{z^*}+\inner{z}{y^*}\\
 \nonumber
 &=\inner{(z^*,z)}{(y,y^*)},
\end{align}
that is,
\[
 \inner{(z^*,z)}{(y,y^*)}-f(y,y^*)\leq \inf_{u^*\in X^*} h^*(u^*,Az)+\delta_{\{0\}}(z^*-A^*u^*).
\]
Now, taking the supremum over $(y,y^*)\in Y\times Y^*$ on the left hand side of the above inequality we obtain
\begin{equation}\label{eq:ap.fs}
 f^*(z^*,z)\leq \inf_{u^*\in X^*} h^*(u^*,Az)+\delta_{\{0\}}(z^*-A^*u^*).
\end{equation}
For a fixed $(z,z^*)\in Y\times Y^*$ such that $f^*(z^*,z)< \infty$, define
$\varphi,\psi\in \BR^{Y\times X\times Y^*\times X^*}$,
\begin{align}
\nonumber
\varphi(y,x,y^*,x^*):=&
 f^*(z^*,z)-\inner{y}{z^*}-\inner{z}{y^*+A^*x^*}+\delta_{\{0\}}(y^*)+h(x,x^*),\\
\nonumber
\psi(y,x,y^*,x^*):=&\delta_{\{0\}}(x-Ay).
\end{align}
Direct calculations yields
\begin{equation}\label{eq:ap.dom}
 \bigcup_{\lambda>0}\lambda [\mbox{D}(\varphi)-\mbox{D}(\psi)]=
  Y\times \bigcup_{\lambda>0}\lambda[\mbox{Pr}_X\mbox{D}(h)]
 -\mbox{R}(A)] \times Y^*\times X^*.
\end{equation}
Using~\eqref{eq:cond2},~\eqref{eq:ap.dom} and Theorem~\ref{th:fr} for $\varphi$ and $\psi$, we conclude that there exists $(y^*,x^*,y^{**},x^{**})\in Y^*\times X^*\times Y^{**}\times X^{**}$
such that
\begin{equation}\label{eq:ap.du}
 \inf \varphi+\psi=-\varphi^*(y^*,x^*,y^{**},x^{**})-
 \psi^*(-y^*,-x^*,-y^{**},-x^{**}).
\end{equation}
Now, notice that
\begin{equation}\label{eq:ap.pos}
 (\varphi+\psi)(y,x,y^*,x^*)\geq
 f^*(z^*,z)+f(y,A^*x^*)-\inner{(z^*,z)}{(y,A^*x^*)}\geq 0.
\end{equation}
Using~\eqref{eq:ap.du} and~\eqref{eq:ap.pos} we get
\begin{equation}\label{eq:ap.leq}
 \varphi^*(y^*,x^*,y^{**},x^{**})+\psi^*(-y^*,-x^*,-y^{**},-x^{**})\leq 0.
\end{equation}
Direct calculations yields
\begin{align}\label{eq:ap.psicj}
 \nonumber
\psi^*(-y^*,-x^*,-y^{**},-x^{**})
=& \sup_{(y,z^*,w^*)}\inner{y}{-y^*-A^*x^*}+\inner{z^*}{-y^{**}}+\inner{w^*}{-x^{**}}\\
=&\, \delta_{\{0\}}(y^*+A^*x^*)+\delta_{\{0\}}(y^{**})+\delta_{\{0\}}(x^{**}).
\end{align}
Now, using~\eqref{eq:ap.leq} and~\eqref{eq:ap.psicj} we conclude that
\[
 y^{**}=0,\,x^{**}=0\quad \mbox{and}\quad y^*=-A^*x^*. 
\]
Therefore, from~\eqref{eq:ap.leq} we have
\begin{align}
\nonumber
\varphi^*(-A^*x^*,x^*,0,0)=&\sup_{(y,x,w^*)}\Big(\inner{y}{z^*-A^*x^*}+\inner{x}{x^*}
+\inner{Az}{w^*}-h(x,w^*)\Big)-f^*(z^*,z)\\
\nonumber
=&\, h^*(x^*,Az)+\delta_{\{0\}}(z^*-A^*x^*)-f^*(z^*,z)\leq 0,
\end{align}
that is, there exists $x^*\in X^*$ such that
\[
 f^*(z^*,z)\geq h^*(x^*,Az)+\delta_{\{0\}}(z^*-A^*x^*).
\]
Finally, using~\eqref{eq:ap.fs} we conclude the proof.
\end{proof}

\end{document}